# SMARANDACHE TYPE FUNCTION OBTAINED BY DUALITY

C. Dumitrescu, N. Vîrlan, Şt. Zamfir, E. Rădescu, N. Rădescu, F.Smarandache
Department of Mathematics, University of Craiova, Romania

**Abstract.** In this paper we extended the Smarandache function from the set $N^*$ of positive integers to the set $Q$ of rational numbers.

Using the inversion formula, this function is also regarded as a generating function. We put in evidence a procedure to construct a (numerical) function starting from a given function in two particular cases. Also connections between the Smarandache function and Euler's totient function as with Riemann's zeta function are established.

## 1. Introduction

The Smarandache function [13] is a numerical function $S : N^* \to N^*$ defined by $S(n) = \min\{m \mid m! \text{ is divisible by } n\}$.

From the definition it results that if

$$n = p_1^{\alpha_1} \cdot p_2^{\alpha_2} \cdots p_t^{\alpha_t} \tag{1}$$

is the decomposition of $n$ into primes, then

$$S(n) = \max S(p_i^{\alpha_i}) \tag{2}$$

and moreover, if $[n_1, n_2]$ is the smallest common multiple of $n_1$ and $n_2$, then

$$S([n_1, n_2]) = \max\{S(n_1), S(n_2)\} \tag{3}$$

The Smarandache function characterizes the prime in the sense that a positive integer $p \geq 4$ is prime if and only if it is a fixed point of $S$.

From Legendre's formula:

$$m! = \prod_p p^{\sum_{i \geq 1} \left[\frac{m}{p^i}\right]} \tag{4}$$

it results [2] that if $a_n(p) = \dfrac{(p^n - 1)}{(p-1)}$ and $b_n(b) = p^n$, then considering the standard numerical scale

$$[p] : b_0(p), b_1(p), \ldots, b_n(p), \ldots$$

$$[p] : a_0(p), a_1(p), \ldots, a_n(p), \ldots$$

we have

$$S(p^k) = p(\alpha_{[p]})_{(p)} \tag{5}$$

that is $S(p^k)$ is calculated multiplying by $p$ the number obtained writing the exponent $\alpha$ in the generalized scale $[p]$ and "reading" it in the standard scale $(p)$.

Let us observe that the calculus in the generalized scale $[p]$ is essentially different from the calculus in the usual scale $(p)$, because the usual relationship $b_{n+1}(p) = pb_n(p)$ is modified in $a_{n+1}(p) = pa_n(p) + 1$ (for more details see [2]).

Let us note from now on $S_p(\alpha) = S(p^\alpha)$. In [3] it is proved that

$$S_p(\alpha) = (p-1)\alpha + \sigma_{[p]}(\alpha) \tag{6}$$

where $\sigma_{[p]}(\alpha)$ is the sum of the digits of $\alpha$ written in the scale $[p]$, and also that

$$S_p(\alpha) = \frac{(p-1)^2}{p}\left(E_p(\alpha) + \alpha\right) + \frac{p-1}{p}\sigma_{(p)}(\alpha) + \sigma_{[p]}(\alpha) \tag{7}$$

where $\sigma_{(p)}(\alpha)$ is the sum of the digits of $\alpha$ written in the standard scale $(p)$ and $E_p(\alpha)$ is the exponent of $p$ in the decomposition into primes of $\alpha!$. From (4) it results that $E_p(\alpha) = \sum_{i \geq 1}\left[\frac{\alpha}{p^i}\right]$, where $[h]$ is the integral part of $h$. It is also said [11] that

$$E_p(\alpha) = \frac{\alpha - \sigma_{(p)}(\alpha)}{p-1} \tag{8}$$

We can observe that this equality may be written as

$$E_p(\alpha) = \left(\left[\frac{\alpha}{p}\right]_{(p)}\right)_{[p]}$$

that is, the exponent of $p$ in the decomposition into primes of $\alpha!$ is obtained writing the integral part of $\alpha/p$ in the base $(p)$ and reading in the scale $[p]$.

Finally, we note that in [1] it is proved that

$$S_p(\alpha) = p\left(\alpha - \left[\frac{\alpha}{p}\right] + \left[\frac{\sigma_{[p]}(\alpha)}{p}\right]\right) \tag{9}$$

From the definition of $S$ it results that $S_p(E_p(\alpha)) = p\left[\frac{\alpha}{p}\right] = \alpha - \alpha_p$ ($\alpha_p$ is the remainder of $\alpha$ with respect to the modulus $m$) and also that

$$E_p(S_p(\alpha)) \geq \alpha \; ; \; E_p(S_p(\alpha) - 1) < \alpha \tag{10}$$

so

$$\frac{S_p(\alpha) - \sigma_{(p)}(S_p(\alpha))}{p-1} \geq \alpha \; ; \; \frac{S_p(\alpha) - 1 - \sigma_{(p)}(S_p(\alpha) - 1)}{p-1} < \alpha.$$

Using (6) we obtain that $S_p(\alpha)$ is the unique solution of the system

$$\sigma_{(p)}(x) \leq \sigma_{[p]}(\alpha) \leq \sigma_{(p)}(x-1) + 1 \tag{11}$$

## 2. Connections with classical numerical functions

It is known that Riemann's zeta function is

$$\zeta(s) = \sum_{n \geq 1} \frac{1}{n^s}.$$

We may establish a connection between the function $S_p$ and Riemann's function as follows:

**Proposition 2.1.** *If* $n = \prod_{i=1}^{t_n} p_i^{\alpha_{i_n}}$ *is the decomposition into primes of the positive integer n then*

$$\frac{\zeta(s-1)}{\zeta(s)} = \sum_{n \geq 1} \prod_{i=1}^{t_n} \frac{S_{p_i}\left(p_i^{\alpha_{i_n}-1}\right) - p_i}{p_i^{s\alpha_{i_n}}}$$

**Proof.** We first establish a connection with Euler's totient function $\varphi$. Let us observe that, for $\alpha \geq 2$, $p^{\alpha-1} = (p-1)a_{\alpha-1}(p) + 1$, so $\sigma_{[p]}(p^{\alpha-1}) = p$. Then by using (6) it results (for $\alpha \geq 2$) that

$$S_p(p^{\alpha-1}) = (p-1)p^{\alpha-1} + \sigma_{[p]}(p^{\alpha-1}) = \varphi(p^\alpha) + p$$

Using the well known relation between $\varphi$ and $\zeta$ given by

$$\frac{\zeta(s-1)}{\zeta(s)} = \sum_{n \geq 1} \frac{\varphi(n)}{n^n}$$

and (12), it results the required relation.

Let us remark also that, if $n$ is given by (1), then

$$\varphi(n) = \prod_{i=1}^{t} \varphi(p_i^{\alpha_i}) = \prod_{i=1}^{t} \left(S_{p_i}(p_i^{\alpha_i-1}) - p_i\right)$$

and

$$S(n) = \max\left(\varphi(p_i^{\alpha_i+1}) + p_i\right)$$

Now it is known that $1 + \varphi(p_i) + \ldots + \varphi(p_i^{\alpha_i}) = p_i^{\alpha_i}$ and then

$$\sum_{k=1}^{\alpha_i-1} Sp_i\left(p_i^k\right) - (\alpha_i - 1)p_i = p_i^{\alpha_i}.$$

Consequently we may write

$$S(n) = \max\left(s \sum_{k=0}^{\alpha_i-1} Sp_i\left(p_i^k\right) - (\alpha_i - 1)p_i\right).$$

To establish a connection with Mangolt's function let us note $\wedge = \min$, $\vee = \max$, $\underset{d}{\wedge} =$ the greatest common divisor, and $\overset{d}{\vee} =$ the smallest common multiple.

We shall write also $n_1 \underset{d}{\wedge} n_2 = (n_1, n_2)$ and $n_1 \overset{d}{\vee} n_2 = [n_1, n_2]$.

The Smarandache function $S$ may be regarded as function from the lattice $\mathcal{L}_d = \left(\mathbb{N}^*, \underset{d}{\wedge}, \overset{d}{\vee}\right)$, into lattice $\mathcal{L} = \left(\mathbb{N}^*, \wedge, \vee\right)$ such that

$$S\left(\underset{i=1,k}{\vee} n_i\right) = \underset{i=1,k}{\vee} S(n_i) \tag{12}$$

Of course $S$ is also order preserving in the sense that $n_1 \leq_d n_2 \to S(n_1) < S(n_2)$.

It is known from [10] that if $(V, \wedge, \vee)$ is a finite lattice, $V = \{x_1, x_2, ..., x_n\}$ with the induced order $\leq$, then for every function $f : V \to \mathbb{N}$ the associated generating function is defined by

$$F(x) = \sum_{y \leq x} f(y) \tag{13}$$

Mangolt's function $\Lambda$ is

$$\Lambda(n) = \begin{cases} \ln p & \text{if } n = p^i \\ 0 & \text{otherwise} \end{cases}$$

The generating function of $\Lambda$ in the lattice $\mathcal{L}_d$ is

$$F^d(n) = \sum_{k \leq_d n} \Lambda(k) = \ln n \tag{14}$$

The last equality follows from the fact that
$$k \leq_d n \Leftrightarrow k \wedge_d n = k \Leftrightarrow k \setminus n \ (k \text{ divides } n)$$

The generating function of $\Lambda$ in the lattice $\mathcal{L}$ is the function $\Psi$

$$F(n) = \sum_{k \leq n} \Lambda(k) = \Psi(n) = \ln[1, 2, ..., n] \tag{15}$$

Then we have the diagram from below.

We observe that the definition of $S$ is in a closed connection with the equalities (1.1) and (2.2) in this diagram. If we note the Mangolt's function by $f$ then the relations

$$[1, 2, ..., n] = e^{F(n)} = e^{f(1)} e^{f(2)} ... e^{f(n)} = e^{\Psi(n)}$$

$$n! = e^{\tilde{F}} = e^{F^d(1)} e^{F^d(2)} ... e^{F^d(n)}$$

together with the definition of $S$, suggest us to consider numerical functions of the form:

$$\nu(n) = \min\{m / n \leq_d [1, 2, ..., m]\} \tag{16}$$

which will be detailed in section 5.

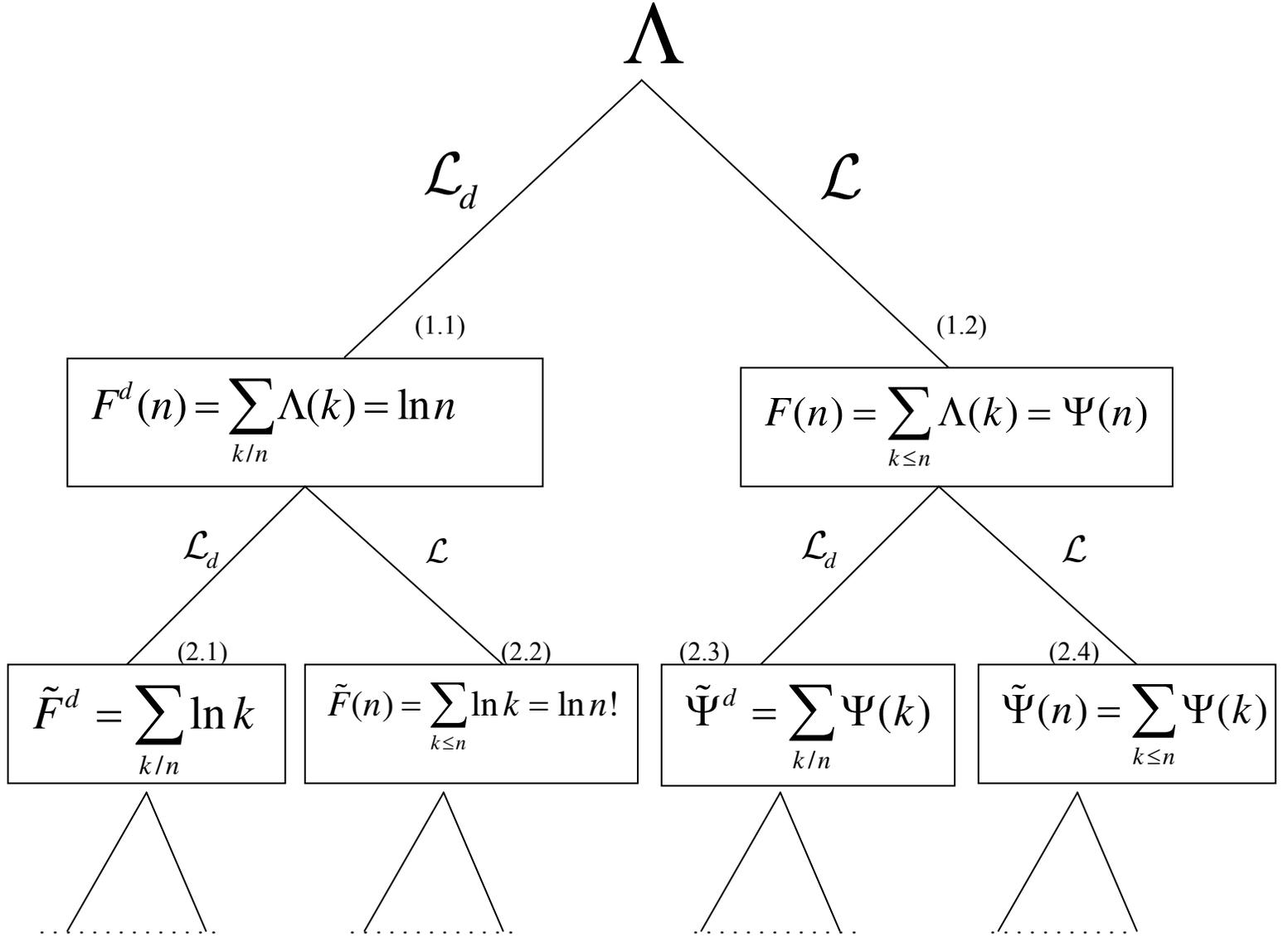

## 3. The Smarandache function as generating function

Let $V$ be a partial order set. A function $f : V \to \mathbb{N}$ may be obtained from its generating function $F$, defined as in (15), by the inversion formula

$$f(x) = \sum_{z \leq x} F(z) \mu(z, x) \qquad (17)$$

where $\mu$ is Moebius function on $V$, that is $\mu : V \times V \to \mathbb{N}$ satisfies:

$$(\mu_1) \, \mu(x, y) = 0, \text{ if } x \not\leq y$$
$$(\mu_2) \, \mu(x, x) = 1$$
$$(\mu_3) \sum_{x \leq y \leq z} \mu(x, y) = 0, \text{ if } x < z.$$

It is known from [10] that if $V = \{1, 2, \ldots, n\}$ then for $(V, \leq_d)$ we have $\mu(x, y) = \mu\left(\dfrac{y}{x}\right)$, where $\mu(k)$ is the numerical Moebius function $\mu(1) = 1$, $\mu(k) = (-1)^i$ if $k = p_1 p_2 \ldots p_k$ and $\mu(k) = 0$ if $k$ is divisible by the square of an integer $d > 1$.

If $f$ is the Smarandache function it results

$$F_S(n) = \sum_{d/n} S(n).$$

Until now it is not known a closed formula for $F_S$, but in [8] it is proved that

(i) $F_S(n) = n$ if and only if $n$ is prime, $n = 9$, $n = 16$, or $n = 24$.

(ii) $F_S(n) > n$ if and only if $n \in \{8, 12, 18, 20\}$ or $n = 2p$ with $p$ a prime (hence it results $F_S(n) \leq n + 4$ for every positive integer $n$) and in [2] it is shown that

$$(iii) \, F(p_1 p_2 \ldots p_t) = \sum_{i=1}^{t} 2^{i-1} p_i.$$

In this section we shall consider the Smarandache function as a generating function, that is using the inversion formula; we shall construct the function $s$ such that

$$s(n) = \sum_{d/n} \mu(d) S\left(\dfrac{n}{d}\right). \qquad (18)$$

If $n$ is given by (1) it results that

$$s(n) = \sum_{p_{i_1} p_{i_2} \ldots p_{i_r}} (-1)^r S\left(\dfrac{n}{p_{i_1} p_{i_2} \ldots p_{i_r}}\right).$$

Let us consider $S(n) = \max S(p_i^{\alpha_i}) = S(p_{i_0}^{\alpha_{i_0}})$. We distinguish the following cases:

($a_1$) if $S(p_{i_0}^{\alpha_{i_0}}) \geq S(p_i^{\alpha_i})$ for all $i \neq i_0$ then we observe that the divisors $d$ for which $\mu(d) \neq 0$ are of the form $d = 1$ or $d = p_{i_1} p_{i_2} \ldots p_{i_r}$. A divisor of the last form may contain $p_{i_0}$ or not, so using (2) it results

$$s(n) = S\left(p_{i_0}^{\alpha_{i_0}}\right)\left(1 - C_{t-1}^1 + C_{t-1}^2 + \ldots + (-1)^{t-1} C_{t-1}^{t-1}\right) + S\left(p_{i_0}^{\alpha_{i_0}-1}\right)\left(-1 + C_{t-1}^1 - C_{t-1}^2 + \ldots + (-1)^t C_{t-1}^{t-1}\right)$$

that is $s(n) = 0$ if $t \geq 2$ or $S\left(p_{i_0}^{\alpha_{i_0}-1}\right)$ and $s(n) = p_{i_0}$ otherwise.

($a_2$)      if there exists $j_0$ such that $S\left(p_{i_0}^{\alpha_{i_0}-1}\right) < S\left(p_{j_0}^{\alpha_{j_0}}\right)$ and

$$S\left(p_{j_0}^{\alpha_{j_0}-1}\right) \geq S\left(p_i^{\alpha_i}\right) \text{ for } i \neq i_0, j_0$$

we also suppose that $S\left(p_{j_0}^{\alpha_{j_0}}\right) = \max\left\{S\left(p_j^{\alpha_j}\right) / S\left(p_{i_0}^{\alpha_{i_0}-1}\right) < S\left(p_j^{\alpha_j}\right)\right\}$.

Then
$$s(n) = S\left(p_{i_0}^{\alpha_{i_0}}\right)\left(1 - C_{t-1}^1 + C_{t-1}^2 - \ldots + (-1)^{t-1} C_{t-1}^{t-1}\right) + S\left(p_{j_0}^{\alpha_{j_0}}\right)\left(-1 + C_{t-2}^1 - C_{t-2}^2 - \ldots + (-1)^{t-1} C_{t-2}^{t-2}\right) +$$
$$+ S\left(p_{j_0}^{\alpha_{j_0}-1}\right)\left(1 - C_{t-2}^1 + C_{t-2}^2 - \ldots + (-1)^{t-2} C_{t-2}^{t-2}\right)$$

so $s(n) = 0$ if $t \geq 3$ or $S\left(p_{j_0}^{\alpha_{j_0}-1}\right) = S\left(p_{j_0}^{\alpha_{j_0}}\right)$ and $s(n) = -p_{j_0}$ otherwise.

Consequently, to obtain $s(n)$ we construct as above a maximal sequence $i_1, i_2, \ldots, i_k$, such that $S(n) = S\left(p_{i_1}^{\alpha_{i_1}}\right)$, $S\left(p_{i_1}^{\alpha_{i_1}-1}\right) < S\left(p_{i_2}^{\alpha_{i_2}}\right), \ldots, S\left(p_{i_{k-1}}^{\alpha_{i_{k-1}}-1}\right) < S\left(p_{i_k}^{\alpha_{i_k}}\right)$ and it results that $s(n) = 0$ if $t \geq k+1$ or $S\left(p_{i_k}^{\alpha_{i_k}}\right) = S\left(p_{i_k}^{\alpha_{i_k}-1}\right)$ and $s(n) = (-1)^{k+1}$ otherwise.

Let us observe that
$$S(p^\alpha) = S(p^{\alpha-1}) \Leftrightarrow (p-1)\alpha + \sigma_{[p]}(\alpha) = (p-1)(\alpha-1) + \sigma_{[p]}(\alpha-1) \Leftrightarrow \sigma_{[p]}(\alpha-1) - \sigma_{[p]}(\alpha) = p - 1$$
Otherwise we have $\sigma_{[p]}(\alpha-1) - \sigma_{[p]}(\alpha) = -1$. So we may write
$$s(n) = \begin{cases} 0 \text{ if } t \geq k+1 \text{ or } \sigma_{[p]}(\alpha_k - 1) - \sigma_{[p]}(\alpha_k) = p - 1 \\ (-1)^{k+1} p_k \text{ otherwise} \end{cases}$$

**Application.** It is known from [10] that $(V, \wedge, \vee)$ is a finite lattice, with the induced order $\leq$ and for the function $f : V \to \mathbb{N}$ we consider the generating function $F$ defined as in (15) then if $g_{ij} = F(x_i \wedge x_j)$ it results $\det g_{ij} = f(x_1) \cdot f(x_2) \cdot \ldots \cdot f(x_n)$. In [10] it is shown also that this assertion may be generalized for partial ordered set by defining
$$g_{ij} = \sum_{\substack{x \leq x_i \\ x \leq x_j}} f(x).$$

Using these results, if we denote by $(i, j)$ the greatest common divisor of $i$ and $j$, and $\Delta(r) = \det\left(S((i,j))\right)$ for $i, j = \overline{1,r}$ then $\Delta(r) = s(1) \cdot s(2) \cdot \ldots \cdot s(r)$. That is for a sufficient large $r$ we have $\Delta(r) = 0$ (in fact for $r \geq 8$). Moreover, for every $n$ there exists a sufficient large $r$ such that $\Delta(n,r) = \det\left(S(n+i, n+j)\right) = 0$, for $i, j = \overline{1,r}$ because $\Delta(n,r) = \prod_{i=1}^{n} S(n+1)$.

## 4. The extension of S to the rational numbers

To obtain this extension we shall define first a dual function of the Smarandache function.

In [4] and [6] a duality principle is used to obtain, starting from a given lattice on the unit interval, other lattices on the same set. The results are used to propose a definition of bi-topological spaces and to introduce a new point of view for studying the fuzzy sets. In [5] the method to obtain new lattices on the unit interval is generalized for an arbitrary lattice.

Here we adopt a method from [5] to construct all the functions tied in a certain sense by duality to the Smarandache function.

Le us observe that if we note $\mathfrak{R}_d(n) = \{m / n \leq_d m!\}$, $\mathcal{L}_d(n) = \{m / m! \leq_d n\}$, $\mathfrak{R}(n) = \{m / n \leq m!\}$, $\mathcal{L}(n) = \{m / m! \leq n\}$ then we may say that the function $S$ is defined by the triplet $(\wedge, \in, \mathfrak{R}_d)$, because $S(n) = \wedge\{m / m \in \mathfrak{R}_d(n)\}$. Now we may investigate all the functions defined by means of a triplet $(a,b,c)$, where $a$ is one of the symbols $\vee, \wedge, \overset{d}{\wedge}, \underset{d}{\vee}$, $b$ is one of the symbols $\in$ and $\notin$, and $c$ is one of the sets $\mathfrak{R}_d(n), \mathcal{L}_d(n), \mathfrak{R}(n), \mathcal{L}(n)$ defined above.

Not all of these functions are non-trivial. As we have already seen the triplet $(\wedge, \in, \mathfrak{R}_d)$ defined the function $S_1(n) = S(n)$, but the triplet $(\wedge, \in, \mathcal{L}_d)$ defines the function $S_2(n) = \wedge\{m / m! \leq_d n\}$, which is identically one.

Many of the functions obtained by this method are step functions. For instance let $S_3$ be the function defined by $(\wedge, \in, \mathfrak{R})$. We have $S_3(n) = \wedge\{m / n \leq m!\}$, so $S_3(n) = m$ if and only if $n \in [(m-1)!+1, m!]$. Let us focus the attention on the function defined by $(\wedge, \in, \mathcal{L}_d)$

$$S_4(4) = \vee\{m / m! \leq_d n\} \tag{19}$$

where there is, in a certain sense, the dual of Smarandache function.

**Proposition 4.1.** *The function $S_4$ satisfies*

$$S_4\left(n_1 \underset{d}{=} \vee n_2\right) = S_4(n_1) \vee S_4(n_2) \tag{20}$$

*so it is a morphism from $\left(\mathbf{N}^*, \underset{d}{\vee}\right)$ to $\left(\mathbf{N}^*, \vee\right)$.*

**Proof.** Let us denote by $p_1, p_2, ..., p_i, ...$ the sequence of the prime numbers and let
$$n_1 = \prod p_i^{\alpha_i}, \quad n_2 = \prod p_i^{\beta_i}.$$

The $n_1 \underset{d}{\wedge} n_2 = \prod p_i^{\min(\alpha_i, \beta_i)}$. If $S_4\left(n_1 \underset{d}{\vee} n_2\right) = m$, $S_4(n_i) = m_i$, for $i = 1, 2$ and we suppose $m_1 \leq m_2$ then the right hand in (22) is $m_1 \wedge m_2 = m$. By the definition $S_4$ we have $E_{p_i}(m) \leq \min(\alpha_i, \beta_i)$ for $i \geq 1$ and there exists $j$ such that $E_{p_i}(m+1) > \min(\alpha_i, \beta_i)$. Then $\alpha_i > E_{p_i}(m)$ and $\beta_i \geq E_{p_i}(m)$ for all $i \geq 1$. We also have $E_{p_i}(m_r) \leq \alpha_i$ for $r = 1,2$. In addition there exist $h$ and $k$ such that $E_{p_h}(m+1) > \alpha_h$, $e_{p_j}(m+1) > \alpha_k$.

Then $\min(\alpha_i, \beta_i) \geq \min(\varepsilon_{p_i}(m_1), \varepsilon_{p_i}(m_2)) = E_{p_i}(m_1)$, because $m_1 \leq m_2$, so $m - 1 \leq m$. If we assume $m_1 < m$ it results that $m! \leq n_1$, therefore it exists $h$ for which $E_{p_h}(m) > \alpha_h$ and we have the contradiction $E_{p_h}(m) > \min\{\alpha_h, \beta_h\}$. Of course $S_4(2n+1) = 1$ and

$$S_4(n) > 1 \text{ if and only if } n \text{ is even.} \tag{21}$$

**Proposition 4.2.** *Let $p_1, p_2, ..., p_i, ...$ be the sequence of all consecutive primes and*
$$n = p_1^{\alpha_1} \cdot p_2^{\alpha_2} \cdot ... \cdot p_k^{\alpha_k} \cdot q_1^{\beta_1} \cdot q_2^{\beta_2} \cdot ... \cdot q_r^{\beta_r}$$
*the decomposition of $n \in \mathbb{N}^*$ into primes such that the first part of the decomposition contains the (eventually) consecutive primes, and let*

$$t_i = \begin{cases} S(p_i^{\alpha_i}) - 1 & \text{if } E_{p_i}(S(p_i^{\alpha_i})) > \alpha_i \\ S(p_i^{\alpha_i}) + p_i - 1 & \text{if } E_{p_i}(S(p_i^{\alpha_i})) = \alpha_i \end{cases} \tag{22}$$

*then $S(n) = \min\{t_1, t_2, ..., t_k, p_{k+1} - 1\}$.*

**Proof.** If $E_{p_i}(S(p_i^{\alpha_i})) > \alpha_i$, then from the definition of the function $S$ results that $S(p_i^{\alpha_i}) - 1$ is the greatest positive integer $m$ such that $E_{p_i}(m) \leq \alpha_i$. Also if $E_{p_i}(S(p_i^{\alpha_i})) = \alpha_i$ then $S(p_i^{\alpha_i}) + p_i - 1$ is the greatest integer $m$ with the property that $E_{p_i}(m) = \alpha_i$.

It results that $\min\{t_1, t_2, ..., t_k, p_{k+1} - 1\}$ is the greatest integer $m$ such that $E_{p-i}(m!) \leq \alpha_i$, for $i = 1, 2, ..., k$.

**Proposition 4.3.** *The function $S_4$ satisfies*
$$S_4((n_1 + n_2)) \wedge S_4([n_1, n_2]) = S_4(n_1) \wedge S_4(n_2)$$
*for all positive integers $n_1$ and $n_2$.*

**Proof.** The equality results using (22) from the fact that
$(n_1 + n_2, [n_1, n_2]) = ((n_1, n_2))$.

We point out now some morphism properties of the functions defined by a triplet $(a, b, c)$ as above.

**Proposition 4.4.**

(i) The function $S_5 : \mathbb{N}^* \to \mathbb{N}^*$, $S_5(n) = \overset{d}{\vee}\{m / m! \leq_d n\}$ satisfies
$$S_5\left(n_1 \underset{d}{\wedge} n_2\right) = S_5(n_1) \underset{d}{\wedge} S_5(n_2) = S_5(n_1) \wedge S_5(n_2) \tag{23}$$

(ii) The function $S_6: \mathbb{N}^* \to \mathbb{N}^*$, $S_6(n) = \overset{d}{\vee}\{m / n \leq_d m!\}$ satisfies

$$S_6\left(n_1 \overset{d}{\vee} n_2\right) = S_6(n_1) \overset{d}{\vee} S_6(n_2) \qquad (24)$$

(iii) The function $S_7: \mathbb{N}^* \to \mathbb{N}^*$, $S_7(n) = \overset{d}{\vee}\{m / m! \leq n\}$ satisfies

$$S_7(n_1 \wedge n_2) = S_7(n_1) \wedge S_7(n_2), \quad S_7(n_1 \vee n_2) = S_7(n_1) \vee S_7(n_2). \qquad (25)$$

**Proof.**

(i) Let $A = \{a_i / a_i! \leq_d n_1\}$, $B = \{b_j / b_j! \leq_d n_2\}$, and $C = \{c_k / c_k! \leq_d n_1 \overset{d}{\vee} n_2\}$. Then we have $A \subset B$ or $B \subset A$. Indeed, let $A = \{a_1, a_2, ..., a_h\}$, $B = \{b_1, b_2, ..., b_r\}$ such that $a_i < a_{i+1}$ and $b_j < b_{j+1}$. Then if $a_h < b_r$ it results that $a_i \leq b_r$ for $i = \overline{1, h}$ so $a_i! \leq_d b_r! \leq_d n_2$. That means $A \subset B$. Analogously, if $b_r \leq a_h$ it results $B \subset A$. Of course we have $C = A \cup B$ so if $A \subset B$ it results

$$S_5\left(n_1 \overset{d}{\wedge} n_2\right) = \overset{d}{\vee} c_k = \overset{d}{\vee} a_i = S_5(n_1) = \min\{S_5(n_1), S_5(n_2)\} = S_5(n_1) \wedge S_5(n_2).$$

From (25) it results that $S_5$ is order preserving in $\mathcal{L}_d$ (but not in $\mathcal{L}$, because $m! < m!+1$ but $S_5(m!) = [1, 2, ..., m]$ and $S_5(m!+1) = 1$, because $m!+1$ is odd).

(ii) Let us observe that $S_6(n) = \overset{d}{\vee}\{m / \exists i \in \overline{1, t}$ such that $E_{p_i}(m) < \alpha_i\}$. If $a = \vee\{m / n \leq_d m!\}$ then $n \leq_d (a+1)!$ and $a+1 = \wedge\{m / n \leq_d m!\} = S(n)$, so $S_6(n) = [1, 2, ..., S(n)-1]$.

Then we have $S_6\left(n_1 \overset{d}{\vee} n_2\right) = \left[1, 2, ..., S\left(n_1 \overset{d}{\vee} n_2\right) - 1\right] = [1, 2, ..., S(n_1) \vee S(n_2) - 1]$

and $S_6(n_1) \overset{d}{\vee} S_6(n_2) = [[1, 2, ..., S_6(n_1)-1], [1, 2, ..., S_6(n_2)-1]] = [1, 2, ..., S_6(n_1) \vee S_6(n_2) - 1]$.

(iii) The relations (27) result from the fact that $S_7(n) = [1, 2, ..., m]$ if and only if $n \in [m!, (m+1)!-1]$.

Now we may extend the Smarandache function to the rational numbers. Every positive rational number $a$ possesses a unique prime decomposition of the form

$$a = \prod_p p^{\alpha_p} \qquad (26)$$

with integer exponents $\alpha_p$, of which only a finite number are nonzero. Multiplication of rational numbers is reduced to addition of their integer exponent system. As a consequence of this reduction questions concerning divisibility of rational numbers are reduced to questions concerning ordering of the corresponding exponent system. That is if $b = \prod_p p^{\beta_p}$ then $b$ divides $a$ if and only if $\beta_p \leq \alpha_p$ for all $p$. The greatest common divisor $d$ and the least common multiple $e$ are given by

$$d = (a, b, ...) = \prod_p p^{\min(\alpha_p, \beta_p, ...)}, \quad e = [a, b, ...] = \prod_p p^{\max(\alpha_p, \beta_p, ...)} \qquad (27)$$

Furthermore, the least common multiple of nonzero numbers (multiplicatively bounded above) is reduced by the rule

$$[a,b,...] = \frac{1}{\left(\frac{1}{a}, \frac{1}{b}, ...\right)} \tag{28}$$

to the greatest common divisor of their reciprocal (multiplicatively bounded below).

Of course we may write every positive rational $a$ under the form $a = n/n_1$, with $n$ and $n_1$ positive integers.

**Definition 4.5.** *The extension $S : Q_+^* \to Q_+^*$ of the Smarandache function is defined by*

$$S\left(\frac{n}{n_1}\right) = \frac{S_1(n)}{S_4(n_1)} \tag{29}$$

A consequence of this definition is that if $n_1$ and $n_2$ are positive integers then

$$S\left(\frac{1}{n_1} \overset{d}{\vee} \frac{1}{n_2}\right) = S\left(\frac{1}{n_1}\right) \vee S\left(\frac{1}{n_2}\right) \tag{30}$$

Indeed

$$S\left(\frac{1}{n_1} \overset{d}{\vee} \frac{1}{n_2}\right) = S\left(\frac{1}{n_1 \underset{d}{\wedge} n_2}\right) = \frac{1}{S_4\left(n_1 \underset{d}{\wedge} n_2\right)} = \frac{1}{S_4(n_1) \wedge S_4(n_2)} = \frac{1}{S_4(n_1)} \vee \frac{1}{S_4(n_2)} = S\left(\frac{1}{n_1}\right) \vee S\left(\frac{1}{n_2}\right)$$

and we can immediately deduce that

$$S\left(\frac{n}{n_1} \overset{d}{\vee} \frac{m}{m_1}\right) = (S(n) \vee S(m)) \cdot \left(S\left(\frac{1}{n_1}\right) \vee S\left(\frac{1}{m_1}\right)\right) \tag{31}$$

It results that function $\overline{S}$ defined by $\overline{S}(a) = \dfrac{1}{S\left(\dfrac{1}{a}\right)}$ satisfies

$$\overline{S}\left(n_1 \underset{d}{\wedge} n_2\right) = \overline{S}(n_1) \wedge \overline{S}(n_2) \text{ and}$$

$$\overline{S}\left(\frac{1}{n_1} \underset{d}{\wedge} \frac{1}{n_2}\right) = \overline{S}\left(\frac{1}{n_1}\right) \wedge \overline{S}\left(\frac{1}{n_2}\right) \tag{32}$$

for every positive integers $n_1$ and $n_2$. Moreover, it results that

$$\overline{S}\left(\frac{n_1}{m_1} \underset{d}{\wedge} \frac{n_2}{m_2}\right) = \left(\overline{S}(n_1) \wedge \overline{S}(n_2)\right) \cdot \left(\overline{S}\left(\frac{1}{m_1}\right) \wedge \overline{S}\left(\frac{1}{m_2}\right)\right)$$

and of course the restriction of $\overline{S}$ to the positive integers is $S_4$. The extension of $S$ to all the rationales is given by $S(-a) = S(a)$.

## 5. Numerical functions inspired from the definition of the Smarandache function

We shall use now the equality (21) and the relation (18) to consider numerical functions as the Smarandache function.

We may say that $m!$ is the product of all positive "smaller" than $m$ in the lattice $\mathcal{L}$. Analogously the product $p_m$ of all the divisors of $m$ is the product of all the elements "smaller" than $m$ in the lattice $\mathcal{L}$. So we may consider functions of the form

$$\Theta(n) = \wedge \{m / n \geq_d p(m)\}. \tag{33}$$

It is known that if $m = p_1^{x_1} \cdot p_2^{x_2} \cdot \ldots \cdot p_t^{x_t}$ then the product of all the divisors of $m$ is $p(m) = \sqrt{m^{\tau(m)}}$ where $\tau(m) = (x_1 + 1)(x_2 + 1)\ldots(x_t + 1)$ is the number of all the divisors of $m$.

If $n$ is given as in (1) then $n \geq_d p(m)$ if and only if

$$\begin{aligned}
g_1 &= x_1(x_1 + 1)(x_2 + 1)\ldots(x_t + 1) - 2\alpha_1 \geq 0 \\
g_2 &= x_2(x_1 + 1)(x_2 + 1)\ldots(x_t + 1) - 2\alpha_2 \geq 0 \\
g_t &= x_t(x_1 + 1)(x_2 + 1)\ldots(x_t + 1) - 2\alpha_t \geq 0
\end{aligned} \tag{34}$$

so $\Theta(n)$ may be obtained solving the problem of non linear programming

$$(\min) f = p_1^{x_1} \cdot p_2^{x_2} \cdot \ldots \cdot p_t^{x_t} \tag{35}$$

under the restrictions (37).

The solution of this problem may be obtained applying the algorithm SUMT (Sequential Unconstrained Minimization Techniques) due to Fiacco and Mc Cormick [7].

**Examples**

1. For $n = 3^4 \cdot 5^{12}$, (37) and (38) become $(\min) f(x) = 3^{x_1} 5^{x_2}$ with $x_1(x_1 + 1)(x_2 + 1) \geq 8$, $x_2(x_1 + 1)(x_2 + 1) \geq 24$. Considering the function

$$U(x, n) = f(x) - r \sum_{i=1}^{k} \ln g_i(x), \text{ and the system}$$

$$\sigma U / \sigma x_1 = 0, \quad \sigma U / \sigma x_2 = 0 \tag{36}$$

in [7] it is shown that if the solution $x_1(r)$, $x_2(r)$ cannot be explained from the system we can make $r \to 0$. Then the system becomes $x_1(x_1 + 1)(x_2 + 1) = 8$, $x_2(x_1 + 1)(x_2 + 1) = 24$ with the (real) solution $x_1 = 1$, $x_2 = 3$.

So we have $\min\{m / 3^4 \cdot 5^{12} \leq \rho(m)\} m_0 = 3 \cdot 5^3$.

Indeed $\rho(m_0) = m_0^{\tau(m_0)/2} = m_0^4 = 3^4 \cdot 5^{12} = n$.

2. For $n = 3^2 \cdot 5^{67}$, from the system (39) it results for $x_2$ the equation $2x_2^3 + 9x_2^2 + 7x_2 - 98 = 0$, with the real solution $x_2 \in (2, 3)$. It results $x_1 \in (4/6, 5/7)$. Considering $x_1 = 1$, we observe that for $x_2 = 2$ the pair $(x_1, x_2)$ is not an admissible solution of the problem, but $x_2 = 3$ gives $\Theta(3^2 \cdot 5^7) = 3^4 \cdot 5^{12}$.

3. Generally, for $n = p_1^{\alpha_1} \cdot p_2^{\alpha_2}$, from the system (39) it results the equation
$$\alpha_1 x_2^3 + (\alpha_1 + \alpha_2) \cdot x_2^2 + \alpha_2 x_2 - 2\alpha_2^2 = 0$$
with solutions given by Cartan's formula.

Of course, using "the method of the triplets", as for the Smarandache function, many other functions may be associated to $\Theta$.

For the function $\nu$ given by (18) it is also possible to generate a class of function by means of such triplets.

In the sequel we'll focus the attention on the analogous of the Smarandache function and on its dual in this case.

**Proposition 5.1.** *If $n$ has the decomposition into primes given by (1) then*

(i) $\nu(n) = \max\limits_{i=\overline{1,t}} p_i^{\alpha_i}$

(ii) $\nu\left(n_1 \stackrel{d}{\vee} n_2\right) = \nu(n_1) \vee \nu(n_2)$

**Proof.**

(i) Let $\max p_i^{\alpha_i} = p_u^{\alpha_u}$. Then $p_i^{\alpha_i} \leq p_u^{\alpha_u}$ for all $\overline{1,t}$, so $p_i^{\alpha_i} \leq_d [1,2,...,p_u^{\alpha_u}]$. But $\left(p_i^{\alpha_i}, p_j^{\alpha_j}\right) = 1$ for $i \neq j$ and then $n \leq_d [1,2,...,p_u^{\alpha_u}]$.

Now if for some $m < p_u^{\alpha_u}$ we have $n \leq_d [1,2,...,m]$, it results the contradiction $p_u^{\alpha_u} \leq_d [1,2,...,m]$.

(ii) If $n_1 = \prod p^{\alpha_p}$, $n_2 = \prod p^{\beta_p}$ then $n_1 \stackrel{d}{\vee} n_2 = \prod p^{\max(\alpha_p, \beta_p)}$ so
$$\nu\left(n_1 \stackrel{d}{\vee} n_2\right) = \max p^{\max(\alpha_p, \beta_p)} = \max\left(\max p^{\alpha_p}, \max p^{\beta_p}\right).$$

The function $\nu_1 = \nu$ is defined by means of the triplet $\left(\vee, \in, \mathfrak{R}_{[d]}\right)$, where $\mathfrak{R}_{[d]} = \{m / n \leq_d [1,2,...,m]\}$. Its dual, in the sense of the above section, is the function defined by the triplet $\left(\vee, \in, \mathcal{L}_{[d]}\right)$. Let us note $\nu_4$ this function
$$\nu_4(n) = \vee\{m \mid [1,2,...,m] \leq_d n\}.$$

That is $\nu_4(n)$ is the greatest natural number with the property that all $m \leq \nu_4(n)$ divide n.

Let us observe that a necessary and sufficient condition to have $\nu_4(n) > 1$ is to exist $m > 1$ such that every prime $p \leq m$ divides $n$. From the definition of $\nu_4$ it also results that $\nu_4(n) = m$ if and only if $n$ is divisible by every $i \leq n$ and not by $m + 1$.

**Proposition 5.2.** *The function $\nu_4$ satisfies*
$$\nu_4\left(n_1 \stackrel{d}{\vee} n_2\right) = \nu_4(n_1) \wedge \nu_4(n_2)$$

**Proof.** Let us note $n = n_1 \overset{d}{\wedge} n_2$, $v_4(n) = m$, $v_4(n_i) = m_i$ for $i = 1,2$. If $m_1 = m_1 \wedge m_2$ then we prove that $m = m_1$. From the definition of $v_4$ it results

$$v_4(n_i) = m_i \leftrightarrow [\forall i \leq m_i \to n \text{ is divisible by } i \text{ but not by } m+1]$$

If $m < m_1$ then $m + 1 \leq m_1 \leq m$ so $m + 1$ divides $n_1$ and $n_2$. That is $m + 1$ divides $n$.
If $m > m_1$ then $m_1 + 1 \leq n$, so $m_1 + 1$ divides $n$. But $n$ divides $n_1$, so $m_1 + 1$ divides $n_1$.
If $t_0 = \max\{i \mid j \leq i \Rightarrow n \text{ divides } n\}$ then $v_4(n)$ may be obtained solving the integer programming problem

$$(\max) f = \sum_{i=1}^{t_0} x_i \ln p$$
$$x_i \leq \alpha_i \text{ for } i = \overline{1, t_0} \qquad (37)$$
$$\sum_{i=1}^{t_0} x_i \ln p_i \leq \ln p_{t_0+1}.$$

If $f_0$ is the maximal value of $f$ for above problem, then $v_4(n) = e^{f_0}$.
For instance $v_4(2^3 \cdot 3^2 \cdot 5 \cdot 11) = 6$.
Of course, the function $v$ may be extended to the rational numbers in the same way as Smarandache function.

**REFERENCES**

[1] M. Andrei, I. Bălăcenoiu, C. Dumitrescu, E. Rădescu, V. Șeleacu – A Linear Combination with the Smarandache Function to otain the Identity – Proceedings of The 26[th] Annual Iranian Mathematics Conference, (1995) pp. 437-439.

[2] M. Andrei, C. Dumitrescu, V. Seleacu, L. Tuțescu, Șt. Zamfir – Some Remarks on the Smarandache Function – Smarandache Function J. 4-5 (1994), pp. 1-5.

[3] M. Andrei, C. Dumitrescu, V. Seleacu, L. Tuțescu, Șt. Zamfir – La function de Smarandache, une nouvelle function dans la theorie des nombres - Congress International H. Poincaré, 14-18 May 1994, Nancy, France.

[4] C. Dumitrescu – Treillius sur des ensembles flous. Applications a des espaces topologiques flous – Rev. Roum. Math. Pures Appl. 31, 1986, pp. 667-675.

[5] C. Dumitrescu – Treillis duals. Applications aux ensembles flous – Math. Rev. d'Anal Numer. Et Theor. de l'Approx., 15, 1986, pp. 111-116.

[6] C. Dumitrescu – Dual Structures in the Fuzzy Sets Theory and in the Groups Theory – Itinerant Sem. on Functional Equations Approx. and Convexity, Cluj-Napoca, Romania, 1989, pp. 23-40.

[7] Fiacco and Mc Cormick – Nonlinear Programming. Sequential unconstrained Minimization Technique - New York, J. Wiley, 1968.

[8] P. Gronas – The Solution of the Diophantine equation $\sigma\eta(n) = n$, Smarandache Function J., V. 4-5, No. 1 (1994), pp. 14-16.

[9] H. Hasse – Number Theory – Akademie-Verlag, Berlin, 1979.
[10] L. Lovasz – Combinatorial Problems and Exercises – Akad. Kiado, Budapest, 1979.
[11] P. Radovici-Marculescu – Probleme de teoria elementară a numerelor – Ed. Tehnică, Bucharest, 1986.
[12] E. Rădescu, N. Rădescu, C. Dumitrescu – On the Sumatory Function associated to Smarandache Function, Smarandache Function J., V. 4-5 (1994), pp.17-21.
[13] F. Smarandache – A function in the Number Theory – An. Univ. Timişoara, Ser. Şt. Math. 28 (1980), pp. 79-88.